\documentclass[12pt]{amsart}
\newtheorem{theorem}{Theorem}
\newtheorem{lemma}{Lemma}
\newtheorem{proposition}{Proposition}

\begin{document}
\centerline{\Large The Lichnerowicz theorem on CR manifolds}
\vskip 0.5cm \centerline{\large Elisabetta Barletta\footnote{The
Author acknowledges support from INdAM within the
interdisciplinary project {\em Nonlinear subelliptic equations of
variational origin in contact geometry}. \newline 2000 Mathematics
Subject Classification: 32V20, 53C17, 53C25, 58C40.\newline Key
words: CR manifold, sublaplacian, Bochner formula, pseudohermitian
Ricci tensor.}}
\title{}
\author{}
\maketitle
\begin{abstract} For any compact strictly pseudoconvex
CR manifold $M$ endowed with a contact form $\theta$ we obtain the
Bochner type formula $\frac{1}{2} \Delta_b (|\nabla^H f|^2 ) =
|\pi_H \nabla^2 f |^2+ (\nabla^H f)(\Delta_b f) + \rho (\nabla^H f
, \nabla^H f) + 2 L f$ (involving the sublaplacian $\Delta_b$ and
the pseudohermitian Ricci curvature $\rho$). When $M$ is compact
of CR dimension $n$ and $\rho (X,X) + 2 A(X, J X) \geq k \,
G_\theta (X,X)$, $X \in H(M)$, we derive the estimate $-\lambda
\geq 2nk/(2n-1)$ on each nonzero eigenvalue $\lambda$ of
$\Delta_b$ satisfying ${\rm Eigen}(\Delta_b ; \lambda ) \cap {\rm
Ker}(T) \neq (0)$ where $T$ is the characteristic direction of $d
\theta$.
\end{abstract}
\section{Introduction}
By a well known result by A. Lichnerowicz, \cite{kn:Lic}, and M.
Obata, \cite{kn:Oba}, on any $m$-dimensional compact Riemannian
manifold $(M , g)$ with ${\rm Ric} \geq k \; g$ the first
eigenvalue of the Laplacian satisfies the estimate
\begin{equation}
\label{e:1}
\lambda_1 \geq m k/(m-1),
\end{equation}
with equality if and only if $M$ is isometric to the standard
sphere $S^m$. The proof of \eqref{e:1} relies on the Bochner
formula (cf. e.g. \cite{kn:BGM}, p. 131)
\begin{equation}
\label{e:2}
- \frac{1}{2} \, \Delta (|d f|^2 ) = |{\rm Hess}(f)|^2  - (d f \,
, \, d \Delta f) + {\rm Ric}((df)^\sharp , (df)^\sharp ),
\end{equation}
for any $f \in C^\infty (M)$. On the other hand, given a compact
strictly pseudoconvex CR manifold $M$, with any fixed contact form
$\theta$ one may associate a natural second order differential
operator $\Delta_b$ (the {\em sublaplacian}) which is similar in
many respects to the Laplacian of a Riemannian manifold. Indeed,
$\Delta_b$ is hypoelliptic and (by a result of \cite{kn:MeSj}) has
a discrete spectrum
\[ 0 < - \lambda_1 < - \lambda_2 < \cdots < - \lambda_k < \cdots \;\; \uparrow + \infty .
\]
Also $(M , \theta )$ carries a natural linear connection $\nabla$
(the {\em Tanaka-Webster connection}, cf.
\cite{kn:Tan}-\cite{kn:Web}) preserving the Levi form and the
maximally complex distribution, and resembling to both the
Levi-Civita connection and the Chern connection (in Hermitian
geometry). Moreover the Ricci tensor $\rho$ of $\nabla$ is likely
to play the role of the Ricci curvature in Riemannian geometry. To
give an example, by a result of J. M. Lee, \cite{kn:Lee}, if $\rho
(Z, \overline{Z}) > 0$ for any $Z \in T_{1,0}(M)$, $Z \neq 0$,
then the first Kohn-Rossi cohomology group $H^{0,1} (M ,
\overline{\partial}_b )$ vanishes (as a CR counterpart of the
classical result in \cite{kn:Boc}). It is a natural question
whether we may estimate the spectrum of $\Delta_b$ from below,
under appropriate geometric assumptions (on $\rho$). The first
attempt to bring \eqref{e:1} to CR geometry belongs to A.
Greenleaf, \cite{kn:Gre}. His result is that on any compact
strictly pseudoconvex CR manifold $M$, of CR dimension $n \geq 3$,
one has
\begin{equation}
\label{e:3} - \lambda_1 \geq n C/(n+1)
\end{equation}
provided that
\begin{equation}
\label{e:4} R_{\alpha\overline{\beta}} Z^\alpha \overline{Z}^\beta
+ i(A_{\overline{\alpha}\overline{\beta}} \overline{Z}^\alpha
\overline{Z}^\beta - A_{\alpha\beta} Z^\alpha Z^\beta ) \geq 2 C
g_{\alpha\overline{\beta}} Z^\alpha \overline{Z}^\beta \, ,
\end{equation}
for some constant $C
> 0$. Here $R_{\alpha\overline{\beta}} = \rho (T_\alpha , T_{\overline{\beta}})$
is the {\em pseudohermitian Ricci tensor} while $A_{\alpha\beta}$
is the {\em pseudohermitian torsion} (cf. e.g. \cite{kn:Dra}, p.
102) and $\{ T_\alpha : 1 \leq \alpha \leq n \}$ is a local frame
of the CR structure. The proof of \eqref{e:3} relies on the rather
involved Bochner like formula
\begin{equation}
\label{e:5} \Delta_b (|\nabla^{1,0} f |^2 ) = 2 \sum_{\alpha  ,
\beta} (f_{\alpha\overline{\beta}} f_{\overline{\alpha}\beta} +
f_{\alpha\beta} f_{\overline{\alpha}\overline{\beta}}) + 4i
\sum_\alpha (f_{\overline{\alpha}} f_{0\alpha} - f_\alpha
f_{0\overline{\alpha}} ) +
\end{equation}
\[ + 2 \sum_{\alpha , \beta} R_{\alpha\overline{\beta}}
f_{\overline{\alpha}} f_\beta + 2 i n \sum_{\alpha , \beta}
(A_{\alpha\beta} f_{\overline{\alpha}} f_{\overline{\beta}} -
A_{\overline{\alpha}\overline{\beta}} f_\alpha f_\beta ) + \] \[ +
\sum_\alpha \{ f_{\overline{\alpha}} (\Delta_b f )_\alpha +
f_\alpha (\Delta_b f)_{\overline{\alpha}} \} \] where
$\nabla^{1,0} f = f^\alpha T_\alpha$. Cf. also Chapter 9 in
\cite{kn:DrTo}. Recently, a large number of results were obtained
within CR and pseudohermitian geometry, mainly by analogy to
similar findings in Riemannian geometry (cf. e.g. S. Dragomir et
al., \cite{kn:DrUr}-\cite{kn:DrNi}). On this line of thought, one
scope of this paper is to establish the Bochner like
formula\footnote{Under the conventions in the present paper the
sublaplacian of \cite{kn:Gre} is $- \Delta_b$.}
\begin{equation}
\label{e:7}
\frac{1}{2} \, \Delta_b (|\nabla^H f |^2 ) =
\end{equation}
\[ = |\pi_H \nabla^2 f|^2
+ (\nabla^H f)(\Delta_b f) + \rho (\nabla^H f , \nabla^H f) + 2 L
f,
\]
for any $f \in C^\infty (M)$, where the differential operator $L$
given by
\begin{equation}
\label{e:L} L f \equiv (J \nabla^H f)(T f) - (J \nabla_T \nabla^H
f)(f).
\end{equation}
As an application we shall prove
\begin{theorem} $\,$ \par\noindent Let $M$ be a compact strictly
pseudoconvex CR manifold, of CR dimension $n$. Let $\theta$ be a
contact form on $M$ such that the Levi form $G_\theta$ is positive
definite. Let $\lambda$ be a nonzero eigenvalue of the
sublaplacian $\Delta_b$. Suppose that there is a constant $k > 0$
such that {\rm i)}
\begin{equation}
\label{e:8} \rho (X,X) + 2 A(X, J X) \geq k \, G_\theta (X,X),
\;\;\; X \in H(M),
\end{equation}
and {\rm ii)} there is an eigenfunction $f \in {\rm
Eigen}(\Delta_b ; \lambda )$ such that $T(f) = 0$. Then $\lambda$
satisfies the estimate
\begin{equation}
\label{e:10}
- \lambda \geq 2nk/(2n-1).
\end{equation}
\label{t:1}
\end{theorem}
\noindent Another lower bound on $- \lambda_1$ (in terms of the
diameter of $(M , g_\theta )$, where $g_\theta$ is the Webster
metric) was found in \cite{kn:BaDr} (by using estimates of the
horizontal gradient at a point, rather than $L^2$ methods) as an
extension of the work by Z. Jiaqing \& Y. Hongcang,
\cite{kn:JiHo}, in Riemannian geometry. Although under more
restrictive assumptions our estimate \eqref{e:10} is sharper than
\eqref{e:3}. When $(M , \theta )$ is Sasakian (i.e.
$A_{\alpha\beta} = 0$) A. Greenleaf's assumption \eqref{e:4}
coincides with our \eqref{e:8}.
\par
The Bochner type formula \eqref{e:7} (as compared to Greenleaf's
\eqref{e:5}) presents a closer resemblance to \eqref{e:2} in
Riemannian geometry, perhaps enabling one to look for an analogue
to the result by M. Obata, \cite{kn:Oba}, as well. Restated in the
CR category, the problem is whether equality in \eqref{e:10}
implies that $M$ is CR isomorphic to the sphere $S^{2n+1}$. As it
turns out when $M = S^{2n+1}$ the assumptions in our Theorem
\ref{t:2} (see below) are satisfied if and only if $n = 1$. We
conjecture that any strictly pseudoconvex CR manifold $M$ carrying
a contact form $\theta$ satisfying (\ref{e:8}) for some $k > 0$
and such that i) $- 2nk /(2n-1) \in {\rm Spec}(\Delta_b )$, and
ii) ${\rm Eigen}(\Delta_b ; - 2nk/(2n-1)) \cap {\rm Ker}(T) \neq
(0)$, is CR isomorphic to $S^3$.

\vskip 0.1in {\small {\bf Acknowledgement}. The Author is grateful
to the anonymous Referee who suggested an improvement of the
statement of Theorem \ref{t:1} and emphasized on the range of
applicability of Theorem \ref{t:2}.}
\section{A reminder of CR geometry}
Let $(M , T_{1,0}(M))$ be an oriented CR manifold, of CR dimension
$n$. For a review of the main notions of CR and pseudohermitian
geometry one may see \cite{kn:Dra}. Let $H(M) = {\rm Re} \{
T_{1,0}(M) \oplus T_{0,1}(M) \}$ be the maximally complex
distribution and $J(Z + \overline{Z}) = i (Z - \overline{Z})$, $Z
\in T_{1,0}(M)$, its complex structure. Let $\theta$ be a {\em
pseudohermitian structure} on $M$, i.e. $\theta$ is a differential
$1$-form such that ${\rm Ker}(\theta ) = H(M)$. The {\em Levi
form} is given by $G_\theta (X,Y ) = (d \theta )(X , J Y)$, $X,Y
\in H(M)$. The given CR manifold is {\em nondegenerate}
(respectively {\em strictly pseudoconvex}) if $G_\theta$ is
nondegenerate (respectively positive definite). From now on, let
us assume that $M$ is nondegenerate. Then each pseudohermitian
structure $\theta$ is a contact form i.e. $\Psi = \theta \wedge (d
\theta )^n$ is a volume form on $M$. Let $T$ be the {\em
characteristic direction} of $d \theta$ i.e. the unique globally
defined nowhere zero tangent vector field $T$ on $M$ determined by
$\theta (T) = 1$ and $T \, \rfloor \, d \theta = 0$. Let
$g_\theta$ be the {\em Webster metric} i.e.
\[ g_\theta (X,Y) = G_\theta (X,Y), \;\;\; g_\theta (X,T) =
0,\;\;\; g_\theta (T,T) = 1, \] for any $X,Y \in H(M)$.  $(M ,
g_\theta )$ is a semi-Riemannian manifold. If $M$ is strictly
pseudoconvex and $\theta$ is chosen such that $G_\theta$ is
positive definite (note that $G_{-\theta}$ is negative definite)
then $(M , g_\theta )$ is a Riemannian manifold (whose canonical
Riemannian volume form is $c_n \Psi$, where $c_n = 2^{-n}/ n !$).
\par
Let $M$ be a strictly pseudoconvex CR manifold and $\theta$ a
contact form on $M$ such that the Levi form $G_\theta$ is positive
definite. The {\em sublaplacian} is
\[ \Delta_b f = {\rm div}(\nabla^H f), \;\;\; f \in C^2 (M), \]
where ${\rm div}(X)$ is the divergence of the vector field $X$
(with respect to the Riemannian metric $g_\theta$) and $\nabla^H f
= \pi_H \nabla f$ is the {\em horizontal gradient}. Precisely
$\nabla f$ is the ordinary gradient (i.e. $g_\theta (\nabla f , X)
= X(f)$ for any $X \in T(M)$) and $\pi_H : T(M) \to H(M)$ is the
projection associated to the direct sum decomposition $T(M) = H(M)
\oplus {\mathbb R} T$. Let $\nabla$ be the Tanaka-Webster
connection of $(M , \theta )$ i.e. the unique linear connection on
$M$ obeying to i) $H(M)$ is $\nabla$-parallel, ii) $\nabla
g_\theta = 0$, $\nabla J = 0$, iii) the torsion $T_\nabla$ of
$\nabla$ satisfies
\[ T_\nabla (Z, W) = 0, \;\;\; T_\nabla (Z , \overline{W}) = 2 i
G_\theta (Z , \overline{W}) T, \;\;\; Z, W \in T_{1,0}(M), \]
\[ \tau \circ J + J \circ \tau = 0, \]
where $\tau (X) = T_\nabla (T, X)$, $X \in T(M)$. A strictly
pseudoconvex CR manifold $M$ is a Sasakian manifold (in the sense
of \cite{kn:Bla}, p. 73) if and only if $\tau = 0$. Given two CR
manifolds $M$ and $N$ a {\em CR map} is a $C^\infty$ map $f : M
\to N$ such that $(d_x f) T_{1,0}(M)_x \subseteq
T_{1,0}(N)_{f(x)}$ for any $x \in M$. A {\em CR isomorphism} is a
$C^\infty$ diffeomorphism and a CR map. By a recent result of G.
Marinescu et al., \cite{kn:MaYe}, any Sasakian manifold is CR
isomorphic to a real submanifold of ${\mathbb C}^N$, for some $N
\geq 2$, carrying the induced CR structure.
\section{The Bochner formula}
Let $\{ X_1 , \cdots X_{2 n} \}$ be a local orthonormal (i.e.
$G_\theta (X_j , X_k ) = \delta_{jk}$) frame of $H(M)$, defined on
the open subset $U \subseteq M$. Then
\begin{equation}
\label{e:11}
\Delta_b f = \sum_{j=1}^{2n} \{ X_j^2 f - (\nabla_{X_j} X_j ) f \}
\end{equation}
on $U$. Let $x_0 \in M$ be an arbitrary point. As well known
$H(M)$ and $g_\theta$ are parallel with respect to $\nabla$.
Therefore, by parallel displacement of a given orthonormal frame
$\{ v_1 , \cdots , v_{2n} \} \subset H(M)_{x_0}$ with $v_{\alpha +
n} = J_x v_\alpha$, $1 \leq \alpha \leq n$, along the geodesics of
$\nabla$ issuing at $x_0$ we may build a local orthonormal frame
$\{ X_j \}$ of $H(M)$, defined on an open neighborhood of $x_0$,
such that
\begin{equation}
\label{e:12}
(\nabla_{X_j} X_k )(x_0 ) = 0, \;\;\; 1 \leq j , k \leq 2n.
\end{equation}
Also $X_{\alpha + n} = J X_\alpha$ (as a consequence of $\nabla J
= 0$). Then (by \eqref{e:11}
and $\nabla g_\theta = 0$)
\[ \Delta_b \left( |\nabla^H f|^2 \right) (x_0 ) = \sum_j X_j^2
\left( |\nabla^H f |^2 \right) (x_0 ) = \]
\[ = 2 \sum_j X_j (g_\theta (\nabla_{X_j} \nabla^H f \, , \,
\nabla^H f ))_{x_0} = \]
\[ = 2 \sum_j \{ g_\theta (\nabla_{X_j} \nabla_{X_j} \nabla^H f \,
, \, \nabla^H f) + g_\theta (\nabla_{X_j} \nabla^H f \, , \,
\nabla_{X_j} \nabla^H f ) \}_{x_0} .  \] As $\{ X_j \}$ is
orthonormal, the first term in the above sum is
\[ \sum_{j,k} g_\theta (\nabla_{X_j} \nabla_{X_j} \nabla^H f \, ,
\, X_k ) X_k (f). \] Moreover (by \eqref{e:12})
\[ g_\theta (\nabla_{X_j} \nabla_{X_j} \nabla^H f \, ,
\, X_k )_{x_0} = \]
\[ = \{ X_j (g_\theta (\nabla_{X_j} \nabla^H f \, , \, X_k )) -
g_\theta (\nabla_{X_j} \nabla^H f \, , \, \nabla_{X_j} X_k )
\}_{x_0} = \]
\[ = X_j \left( X_j (g_\theta (\nabla^H f , X_k )) - g_\theta (\nabla^H
f , \nabla_{X_j} X_k ) \right)_{x_0} = \]
\[ = X_j \left( X_j X_k f - (\nabla_{X_j} X_k ) f \right)_{x_0} =
X_j \left( (\nabla^2 f)(X_j , X_k ) \right)_{x_0} \] where the
Hessian is defined with respect to the Tanaka-Webster connection
\[ (\nabla^2 f)(X,Y) = (\nabla_X d f) Y = X(Y(f)) - (\nabla_X Y)
f, \;\;\; X,Y \in T(M). \] Unlike the Hessian in Riemannian
geometry $\nabla^2 f$ is never symmetric
\begin{equation}
\label{e:13}
(\nabla^2 f)(X,Y) = (\nabla^2 f)(Y,X) - T_\nabla (X,Y) (f),
\end{equation}
where $T_\nabla$ is the torsion of $\nabla$. On the other hand
$T_\nabla$ is pure (cf. \cite{kn:Dra}, p. 102) hence
\begin{equation}
\label{e:14}
T_\nabla (X,Y) = - 2 \Omega (X,Y) T, \;\;\; X,Y \in H(M).
\end{equation}
Here $\Omega (X,Y) = g_\theta (X, J Y)$ (so that $\Omega = - d
\theta$). Then (by \eqref{e:13}-\eqref{e:14})
\[ g_\theta (\nabla_{X_j} \nabla_{X_j} \nabla^H f \, , \, X_k
)_{x_0} = X_j \left( (\nabla^2 f)(X_j , X_k ) \right)_{x_0} = \]
\[ = X_j \left( (\nabla^2 f)(X_k , X_j ) + 2 \Omega (X_j , X_k ) T
f \right)_{x_0} = \]
\[ = g_\theta (\nabla_{X_j} \nabla_{X_k} \nabla^H f \, , \, X_j
)_{x_0} + 2 \Omega (X_j , X_k )_{x_0} X_j (T f)_{x_0} \] so that
\begin{equation}
\label{e:15}
\frac{1}{2} \, \Delta_b \left( |\nabla^H f |^2 \right) (x_0 ) =
\sum_{j} \left| \nabla_{X_j} \nabla^H f \right|^2_{x_0} +
\end{equation}
\[ + \sum_{j,k} \{ g_\theta (\nabla_{X_j} \nabla_{X_k} \nabla^H f \, , \, X_j
) + 2 \Omega (X_j , X_k ) X_j (T f) \}_{x_0} \; X_k (f)_{x_0} . \]
If $B$ is a bilinear form on $T(M)$ we denote by $\pi_H B$ its
restriction to $H(M)$. The norm of $\pi_H B$ is given by $|\pi_H
B|^2 = \sum_{j,k} B(X_j , X_k )^2$. Then
\[ |\pi_H \nabla^2 f |^2 = \sum_{j,k} (\nabla^2 f)(X_j , X_k )^2 =
\sum_{j,k} \left( X_j X_k f - (\nabla_{X_j} X_k ) f \right)^2 = \]
\[ = \sum_{j,k} g_\theta (\nabla_{X_j} \nabla^H f \, , \, X_k )^2
= \sum_j g_\theta (\nabla_{X_j} \nabla^H f \, , \, \nabla_{X_j}
\nabla^H f) \] so that
\begin{equation}
\label{e:16}
|\pi_H \nabla^2 f |^2 = \sum_j |\nabla_{X_j} \nabla^H f |^2 .
\end{equation}
Next $[X_j , X_k ] = \nabla_{X_j} X_k - \nabla_{X_k} X_j - T_\nabla
(X_j , X_k )$ hence (by applying \eqref{e:12} and \eqref{e:14})
\[ [X_j , X_k ]_{x_0} = 2 \Omega (X_j , X_k )_{x_0} T_{x_0}  \]
and taking into account
\[ \nabla_X \nabla_Y = \nabla_Y \nabla_X + R(X,Y) + \nabla_{[X,Y]}
\]
(where $R$ is the curvature tensor field of $\nabla$) we obtain
\begin{equation}
\label{e:17}
\nabla_{X_j} \nabla_{X_k} \nabla^H f = \nabla_{X_k} \nabla_{X_j}
\nabla^H f +
\end{equation}
\[ + R(X_j , X_k ) \nabla^H f + 2 \Omega (X_j , X_k ) \nabla_T \nabla^H f \]
at $x_0$. Moreover
\[ g_\theta (\nabla_{X_k} \nabla_{X_j} \nabla^H f \, , \, X_j
)_{x_0} = \] \[ = \{ X_k (g_\theta (\nabla_{X_j} \nabla^H f \, ,
\, X_j )) - g_\theta (\nabla_{X_j} \nabla^H f \, , \, \nabla_{X_k}
X_j ) \}_{x_0} = \]
\[ = X_k \left( X_j^2 f - (\nabla_{X_j} X_j ) f \right)_{x_0} \]
that is
\begin{equation}
\label{e:18}
\sum_j g_\theta (\nabla_{X_k} \nabla_{X_j} \nabla^H f \, , \, X_j
)_{x_0} = X_k \left( \Delta_b f \right)_{x_0}.
\end{equation}
Therefore (by \eqref{e:17}-\eqref{e:18})
\[ \sum_{j,k} g_\theta (\nabla_{X_j} \nabla_{X_k} \nabla^H f \,
, \, X_j )_{x_0} X_k (f)_{x_0} = \] \[ = \sum_k \{ X_k (\Delta_b
f) \; X_k f \}_{x_0} + \sum_{j,k} \{ g_\theta (R(X_j , X_k )
\nabla^H f \, , \, X_j ) \; X_k f + \] \[ + 2 \Omega (X_j , X_k )
g_\theta (\nabla_T \nabla^H f \, , \, X_j ) \; X_k f \}_{x_0} = \]
\[ = (\nabla^H f)(\Delta_b f)_{x_0} + \sum_j \{ g_\theta (R(X_j ,
\nabla^H f ) \nabla^H f \, , \, X_j ) + \] \[ + 2 g_\theta (X_j \,
, \, J \nabla^H f ) g_\theta (\nabla_T \nabla^H f \, , \, X_j )
\}_{x_0} = \] \[ = (\nabla^H f)(\Delta_b f)_{x_0} + \rho (\nabla^H
f \, , \, \nabla^H f)_{x_0} + 2 g_\theta (\nabla_T \nabla^H f \, ,
\, J \nabla^H f)_{x_0} \] where $\rho (X,Y) = {\rm trace} \{ Z
\mapsto R(Z, Y)X \}$. Then (by \eqref{e:16}) the identity
\eqref{e:15} becomes
\[ \frac{1}{2} \, \Delta_b \left( |\nabla^H f |^2 \right) = |\pi_H \nabla^2
f |^2 + (\nabla^H f)(\Delta_b f) + \rho (\nabla^H f , \nabla^H f)
+ \]
\[ + 2 g_\theta (\nabla_T \nabla^H f \, , \, J \nabla^H f) + 2
g_\theta (\nabla^H T f \, , \, J \nabla^H f) \] which yields
\eqref{e:7}.
\section{A lower bound on $-\lambda$ for $\lambda \in {\rm Spec}(\Delta_b )$ with
${\rm Eigen}(\Delta_b ; \lambda ) \cap {\rm Ker}(T) \neq (0)$} Let
$M$ be a compact strictly pseudoconvex CR manifold and $\theta$ a
contact form on $M$ with $G_\theta$ positive definite. Let $(u,v)
= \int_M u v \, \Psi$ be the $L^2$ inner product on $M$ and $\| u
\| = (u,u)^{1/2}$ the $L^2$ norm. For any $f \in C^\infty (M)$ let
$f_0 = T(f)$ . We shall need the following two lemmas.
\begin{lemma}
\label{l:1}
\begin{equation}
\label{e:23} {\rm div}(J \nabla^H f) = 2 n f_0 \, .
\end{equation}
\end{lemma}
{\em Proof.} Let $\{ T_\alpha : 1 \leq \alpha \leq n \}$ be a local frame of
$T_{1,0}(M)$, defined on $U \subseteq M$. Then $\nabla^H f =
f^\alpha T_\alpha + f^{\overline{\alpha}} T_{\overline{\alpha}}$
on $U$, where $f^\alpha = g^{\alpha\overline{\beta}}
f_{\overline{\beta}}$, $f_{\overline{\beta}} =
T_{\overline{\beta}} (f)$ and $T_{\overline{\beta}} =
\overline{T}_\beta$, hence
\begin{equation}
\label{e:20}
J \nabla^H f = i (f^\alpha T_\alpha - f^{\overline{\alpha}}
T_{\overline{\alpha}}) .
\end{equation}
We wish to compute the divergence of the vector field
\eqref{e:20}. As $\Psi$ is parallel with respect to $\nabla$
\[ {\rm div}(J \nabla^H f) = {\rm trace} \{ T_A \mapsto
\nabla_{T_A} J \nabla^H f \} \] where $A \in \{ 0 , 1, \cdots , n
, \overline{1} , \cdots , \overline{n} \}$ (with the convention
$T_0 = T$). We set $f_{AB} = (\nabla^2 f)(T_A , T_B )$. Then
\[ \nabla_{T_\beta} J \nabla^H f = i ({f_\beta}^\alpha T_\alpha -
{f_{\beta}}^{\overline{\alpha}} T_{\overline{\alpha}} ) \] where
${f_\beta}^\alpha = g^{\alpha\overline{\gamma}}
f_{\beta\overline{\gamma}}$, etc., so that
\begin{equation}
\label{e:21}
{\rm div}(J \nabla^H f) = i( {f_\alpha}^\alpha -
{f_{\overline{\alpha}}}^{\overline{\alpha}} ) .
\end{equation}
The identities \eqref{e:13}-\eqref{e:14} furnish the commutation
formula $f_{\alpha\overline{\beta}} = f_{\overline{\beta}\alpha} -
2 i g_{\alpha\overline{\beta}} \, f_0$. In particular
\begin{equation}
\label{e:22} 
{f_{\overline{\alpha}}}^{\overline{\alpha}} = {f_\alpha}^\alpha + 2
i n f_0
\end{equation}
hence \eqref{e:21} yields \eqref{e:23}. Q.e.d.
\begin{lemma}
\label{l:2}
\begin{equation}
\label{e:24}
\int_M L f \,  \Psi = - 4 n  \| f_0 \|^2  + \int_M  A(\nabla^H f \, , \,
J \nabla^H f ) \Psi \, ,
\end{equation}
Here $A (X,Y) = g_\theta (\tau X , Y)$ is the pseudohermitian
torsion of $(M , \theta )$ and $L$ is given by {\rm (\ref{e:L})}.
\end{lemma}
{\em Proof.} By the very definition of $L f$
\[ \int_M L f \, \Psi = J_1 -J_2 \]
where
\[ J_1 = \int_M (J \nabla^H f)(f_0 ) \Psi \;\; , \;\;\;
J_2 = \int_M (J \nabla_T \nabla^H f ) (f) \, \Psi \, .\] By
Green's lemma and \eqref{e:23}
\[ J_1 = \int_M \{ {\rm div}(f_0 \; J
\nabla^H f ) - f_0 \; {\rm div}(J \nabla^H f) \} \Psi = \]
\[ = - \int_M f_0 \; {\rm div}(J \nabla^H f) \Psi =  - 2n \| f_0 \|^2 \, ,\]
\[J_2 = \int_M {\rm div}( f J \nabla_T \nabla^H f ) \Psi -
\int_M f {\rm div}(J \nabla_T \nabla^H f ) \Psi =\]
\[ = - \int_M f \; {\rm div} (J \nabla_T \nabla^H f) \Psi . \]
Let us compute in local coordinates ${\rm div} (J \nabla_T \nabla^H
f)$. According to the notations used in the proof of lemma
\ref{l:1}, set $f_{ABC}
= (\nabla^3 f)(T_A , T_B , T_C )$ where
\[ (\nabla^3 f)(X,Y,Z) = (\nabla_X \nabla^2 f)(Y,Z) = \]
\[ = X((\nabla^2 f)(Y,Z)) - (\nabla^2 f)(\nabla_X Y , Z) - (\nabla^2
f)(Y , \nabla_X Z), \] for any $X,Y, Z \in T(M)$. Then
\[ J \nabla_T \nabla^H f = i({f_0}^\alpha T_\alpha -
{f_0}^{\overline{\alpha}} T_{\overline{\alpha}} ) \] yields
\[ \nabla_{T_\alpha} (J \nabla_T \nabla^H f) = i (
{f_{\alpha 0}}^\beta T_\beta - {f_{\alpha 0}}^{\overline{\beta}}
T_{\overline{\beta}} ) \] so that (by $\nabla T = 0$)
\begin{equation}
\label{e:25}
{\rm div}(J \nabla_T \nabla^H f ) = {\rm trace} \{ T_A \mapsto
\nabla_{T_A} J \nabla_T \nabla^H f \} = i ({f_{\alpha 0}}^\alpha -
{f_{\overline{\alpha}0}}^{\overline{\alpha}} ),
\end{equation}
where ${f_{\alpha 0}}^\beta = g^{\beta\overline{\gamma}} f_{\alpha
0 \overline{\gamma}}$, etc. We need the third order commutation
formula
\begin{equation}
\label{e:26}
f_{\overline{\beta}0\alpha} = f_{\alpha 0 \overline{\beta}} + 2 i
g_{\alpha\overline{\beta}} f_{00} + A^\gamma_{\overline{\beta}}
f_{\alpha \gamma} - A_\alpha^{\overline{\gamma}}
f_{\overline{\beta}\overline{\gamma}} + A^\gamma_{\overline{\beta}
, \alpha}  f_\gamma - A^{\overline{\gamma}}_{\alpha ,
\overline{\beta}} f_{\overline{\gamma}}
\end{equation}
where $f_{00} = (\nabla^2 f)(T,T) = T^2 (f)$. This follows from
\[ (\nabla^3 f)(X,T, Y) - (\nabla^3 f)(Y,T, X) = \]
\[ = 2 \Omega (X,Y) f_{00} - X(\tau (Y)f) + Y(\tau (X)f) + \tau ([X,Y])f \]
i.e.
\[ (\nabla^3 f)(X,T,Y) = (\nabla^3 f)(Y,T,X) + 2 \Omega (X,Y)
f_{00} + \] \[ + (\nabla^2 f)(Y, \tau (X)) - (\nabla^2 f)(X, \tau
(Y)) - S(X,Y) f \, , \] for any $X,Y \in H(M)$, where $S(X,Y) =
(\nabla_X \tau )Y - (\nabla_Y \tau )X$. Indeed we may set $X =
T_\alpha$ and $Y = T_{\overline{\beta}}$ in the previous identity
and observe that $S(T_\alpha , T_{\overline{\beta}}) =
A^\gamma_{\overline{\beta} , \alpha} T_\gamma -
A^{\overline{\gamma}}_{\alpha , \overline{\beta}}
T_{\overline{\gamma}}$ where $\tau (T_\alpha ) =
A^{\overline{\beta}}_\alpha T_{\overline{\beta}}$ and the covariant
derivatives $A^\gamma_{\overline{\beta}, \alpha}$ are given by
$(\nabla_{T_\alpha} \tau ) T_{\overline{\beta}} =
A^\gamma_{\overline{\beta}, \alpha} T_\gamma$. The identity
\eqref{e:26} leads to
\[ {f_{\overline{\alpha}0}}^{\overline{\alpha}} =
{f_{\alpha 0}}^\alpha + 2 i n f_{00} + A^{\alpha\beta}
f_{\alpha\beta} - A^{\overline{\alpha}\overline{\beta}}
f_{\overline{\alpha}\overline{\beta}} + {A^{\alpha\beta}}_{,
\alpha} f_\beta - {A^{\overline{\alpha}\overline{\beta}}}_{,
\overline{\alpha}} f_{\overline{\beta}} \]
hence \eqref{e:25} becomes
\[ {\rm div}(J \nabla_T \nabla^H f) = 2 n f_{00} - i( A^{\alpha\beta}
f_{\alpha\beta} - A^{\overline{\alpha}\overline{\beta}}
f_{\overline{\alpha}\overline{\beta}} + {A^{\alpha\beta}}_{,
\alpha} f_\beta - {A^{\overline{\alpha}\overline{\beta}}}_{,
\overline{\alpha}} f_{\overline{\beta}}) \; . \]
Therefore
\[ J_2 = - 2n \int f f_{00}
\Psi  
+ i \int f (A^{\alpha\beta} f_{\alpha\beta} -
A^{\overline{\alpha}\overline{\beta}}
f_{\overline{\alpha}\overline{\beta}} + {A^{\alpha\beta}}_{,
\alpha} f_\beta - {A^{\overline{\alpha}\overline{\beta}}}_{,
\overline{\alpha}} f_{\overline{\beta}}) \Psi \]
where
\[ \int_M f f_{00} = \int_M f \; T (f_0 ) \Psi = \int_M \{ T(f f_0 ) - f_0^2
\} \Psi =\]
\[ = \int_M \{ {\rm div}(f f_0 T) - f f_0 \; {\rm div}(T) \} \Psi -
\| f_0\|^2 \]
hence (by ${\rm div}(T) = 0$)
\[ \int_M f f_{00} \Psi = - \| f_0 \|^2 \, . \]
On the other hand ${\rm div}(Z^\alpha T_\alpha ) = {Z^\alpha}_{,
\alpha}$ hence (by Green's lemma)
\[ \int_M f {A^{\alpha\beta}}_{, \alpha} f_\beta = \int_M \{ (f f_\beta
A^{\alpha\beta})_{, \alpha} - A^{\alpha\beta} f_\alpha f_\beta - f
A^{\alpha\beta} f_{\beta , \alpha} \} \Psi= \]
\[ = - \int_M ( A^{\alpha\beta} f_\alpha f_\beta + f A^{\alpha\beta} f_{\beta ,
\alpha} ) \Psi \]
where $f_{\beta , \alpha} = (\nabla_{T_\alpha} df) T_\beta =
f_{\alpha\beta}$. Hence
\[ J_2 = 2 n \| f_0\|^2 + i \int_M  (A^{\overline{\alpha}
\overline{\beta}} f_{\overline{\alpha}} f_{\overline{\beta}} -
A^{\alpha\beta} f_\alpha f_\beta ) \Psi \] and then (by $A(\nabla^H
f , J \nabla^H f) = i (A_{\alpha\beta} f^\alpha f^\beta
- A_{\overline{\alpha}\overline{\beta}} f^{\overline{\alpha}}
f^{\overline{\beta}})$) we may conclude that
\[ J_2 = 2 n \| f_0\|^2 - \int_M A( \nabla^H f , J \nabla^H f )
\Psi \] so Lemma \ref{l:2} is proved. \vskip 0.5cm Let us prove
Theorem \ref{t:1}. Note that $(\nabla^H f)(f) = |\nabla^H f|^2$.
Let $\lambda$ be an eigenvalue of $\Delta_b$ and $f$ an
eigenfunction corresponding to $\lambda$ such that $T(f) = 0$.
Then \eqref{e:7} becomes
\[ \frac{1}{2} \; \Delta_b \left( |\nabla^H f |^2 \right) = |\pi_H
\nabla^2 f |^2 + \lambda |\nabla^H f |^2 + \rho (\nabla^H f ,
\nabla^H f ) + 2 \, L f . \] Let us integrate over $M$ and use
Green's lemma, Lemma \ref{l:2} and the assumptions (i)-(ii) in
Theorem \ref{t:1} to get
\[ 0 = \| \pi_H \nabla^2 f \|^2 + \lambda \| \nabla^H f \|^2 + \int_M
\{ \rho (\nabla^H f , \nabla^H f ) + 2 A(\nabla^H f , J \nabla^H
f) \} \Psi \geq
\]
\[ \geq \| \pi_H \nabla^2 f \|^2 + (\lambda + k) \| \nabla^H f \|^2
\] that is
\begin{equation}
\label{e:27}
0 \geq \| \pi_H \nabla^2 f \|^2 + (\lambda + k) \| \nabla^H f \|^2
\, .
\end{equation}
Once again, as $f$ is an eigenfunction
\[ \| \Delta_b f \|^2 = \int_M | \Delta_b f |^2 \Psi = \lambda
\int_M f \; \Delta_b f\, \Psi = \lambda \int_M f \; {\rm
div}(\nabla^H f) \Psi = \]
\[ = \lambda \int_M \{ {\rm div}(f \; \nabla^H f) - (\nabla^H
f)(f) \} \Psi \] that is
\begin{equation}
\label{e:28}
\| \Delta_b f \|^2 = - \lambda \| \nabla^H f \|^2 .
\end{equation}
Next (with the notations in Section 2) we set
\[ v_j = ((\nabla^2 f)(X_j , X_1 ) , \cdots , (\nabla^2 f)(X_j ,
X_{2n}))\; , \;\;\; 1 \leq j \leq 2n \, , \] so that
\[ |\pi_H \nabla^2 f |^2 = \sum_{j,k} (\nabla^2 f )(X_j , X_k )^2
= \sum_j |v_j |^2 = |w|^2 \]
where $w = (|v_1 |, \cdots , |v_{2n}|)$ (and $|v_j|$, $|w|$ are the
Euclidean norm of $v_j$, $w$). By the Cauchy-Schwarz inequality
\[ |\pi_H \nabla^2 f |^2 = |w|^2 \geq \frac{1}{2n} \left| w \cdot
(1, \cdots , 1) \right|^2 = \]
\[ = \frac{1}{2n} \left( \sum_j |v_j | \right)^2 \geq
\frac{1}{2n} \left( \sum_j |(\nabla^2 f)(X_j , X_j )| \right)^2 \]
hence
\begin{equation}
\label{e:29}
|\pi_H \nabla^2 f |^2 \geq \frac{1}{2n} \; (\Delta_b f)^2 .
\end{equation}
Finally, by \eqref{e:28}-\eqref{e:29} the inequality \eqref{e:27}
becomes (whenever $\| \Delta_b f \| \neq 0$, that is $\lambda \neq
0$)
\[ 0 \geq \left( \frac{1}{2n} - \frac{\lambda + k}{\lambda}
\right) \; \| \Delta_b f \|^2 \] to conclude that $-\lambda \geq
2nk/(2n-1)$. Q.e.d.

\vskip 0.1in We close the section with the following remark on
assumption (ii) in Theorem \ref{t:1}. The problem whether ${\rm
Eigen}(\Delta_b ; \lambda ) \cap {\rm Ker}(T) \neq \emptyset$ is
in general open. Nevertheless if $M = S^{2n+1}$ is the standard
sphere then eigenfunctions $f \in {\rm Eigen}(\Delta_b ; -4(n+1))$
with $T(f) = 0$ may be easily produced (here $-4(n+1)$ is the
second nonzero eigenvalue of the ordinary Laplacian on
$S^{2n+1}$). Indeed let $\Delta$ be the Laplace-Beltrami operator
of $(S^{2n+1} , g_{\theta_0})$. As well known (cf. e.g.
\cite{kn:BGM}) $\Delta v = - \ell (\ell + 2n) v$, where $v$ is the
restriction to $S^{2n+1}$ of a harmonic polynomial $H \in
{\mathcal H}_\ell$ (here ${\mathcal H}_\ell$ is the space of
harmonic, i.e. $\Delta_{{\mathbb R}^{2n+2}} H = 0$, polynomials $H
: {\mathbb R}^{2n+2} \to {\mathbb R}$ which are homogeneous of
degree $\ell$) and the whole spectrum of $\Delta$ on $S^{2n+1}$
may be obtained this way. Note that ${\mathcal H}_2$ consists of
all $H = \sum_{i,j} (a_{ij} x^i x^j + b_{ij} x^i y^j + c_{ij} y^i
y^j )$ with $\sum_i (a_{ii} + c_{ii}) = 0$. For the sphere $(d
\iota ) T = T_0$ where $\iota : S^{2n+1} \to {\mathbb C}^{n+1}$ is
the inclusion while $T_0 = x^j
\partial /\partial y^j - y^j \partial /\partial x^j$ and $(x^j ,
y^j )$ are the natural coordinates on ${\mathbb C}^{n+1} \approx
{\mathbb R}^{2n+2}$, hence
\begin{equation}
\label{e:eigenfunction} {\mathcal H}_2 \cap {\rm Ker}(T_0 ) = \{ H
= \sum_{i,j} a_{ij}(x^i x^j + y^i y^j ) : \sum_{i} a_{ii} = 0 \} .
\end{equation}
Finally, by a formula of A. Greenleaf (cf. {\em op. cit.})
\begin{equation} \label{e:Greenleaf}
\Delta_b = \Delta - T^2
\end{equation}
hence $-4(n+1) \in {\rm Spec}(\Delta_b )$ and $(0) \neq {\rm
Eigen}(\Delta ; - 4(n+1)) \cap {\rm Ker}(T) \subseteq {\rm
Eigen}(\Delta_b ; -4(n+1))$. On the other hand note that
${\mathcal H}_1 \cap {\rm Ker}(T_0 ) = (0)$. So the eigenfunctions
of $\Delta_b$ we consider (cf. (\ref{e:eigenfunction}) above) are
spherical harmonics of degree 2. However $4(n+1)$ is greater equal
than minus the third eigenvalue of $\Delta_b$ (cf. Proposition
\ref{p:2} below). See also our Appendix A for a short proof of
(\ref{e:Greenleaf}).
\section{Consequences of $- 2nk/(2n-1) \in {\rm Spec}(\Delta_b )$.}
Let $M$ be a strictly pseudoconvex CR manifold and $\theta$ a
contact form on $M$ such that $G_\theta$ is positive definite. We
recall a few concepts from {\em sub-Riemannian geometry} (cf. e.g.
R. S. Strichartz, \cite{kn:Str}) on a strictly pseudoconvex CR
manifold. Let $x \in M$ and $g(x) : T^*_x (M) \to H(M)_x$
determined by
\[ G_{\theta , x} (v , g(x) \xi ) = \xi (v), \;\;\; v \in H(M)_x ,
\;\; \xi \in T_x^* (M). \] Note that the kernel of $g$ is
precisely the conormal bundle
\[ H(M)^\bot_x = \{ \omega \in T^*_x (M) : {\rm Ker}(\omega ) \supseteq
H(M)_x \} , \;\;\; x \in M. \] That is $G_\theta$ is a {\em
sub-Riemannian metric} on $H(M)$ and $g$ its alternative
description (cf. (2.1) in \cite{kn:Str}, p. 225).
\par
Let $\gamma : I \to M$ be a piecewise $C^1$ curve (where $I
\subseteq {\mathbb R}$ is an interval). Then $\gamma$ is a {\em
lengthy curve} if $\dot{\gamma}(t) \in H(M)_{\gamma (t)}$ for
every $t \in I$ such that $\dot{\gamma}(t)$ is defined. A
piecewise $C^1$ curve $\xi : I \to T^* (M)$ is a {\em cotangent
lift} of $\gamma$ if $\xi (t) \in T_{\gamma (t)}^* (M)$ and
$g(\gamma (t)) \xi (t) = \dot{\gamma}(t)$ for every $t$ (where
defined). The {\em length} of a lengthy curve $\gamma : I \to M$
in sub-Riemannian geometry
\[ L(\gamma ) = \int_I \{ \xi (t) \left[ g(\gamma (t))
\xi (t) \right]\}^{1/2} \; d t = \int_I G_{\theta , \gamma (t)}
(\dot{\gamma}(t) , \dot{\gamma} (t))^{1/2}  \] coincides with the
Riemannian length of $\gamma$ as a curve in $(M , g_\theta )$. The
{\em Carnot-Carath\'eodory distance} $\rho (x,y)$ among $x, y \in
M$ is the infimum of the lengths of all lengthy curves joining $x$
and $y$. By a well known theorem of W. L. Chow, \cite{kn:Cho}, any
two points $x , y \in M$ may be joined by a lengthy curve (and one
may easily check that $\rho$ is a distance function on $M$).
\par
Let $g_\theta$ be the Webster metric of $(M , \theta )$. Then
$g_\theta$ is a {\em contraction} of the sub-Riemannian metric
$G_\theta$ ($G_\theta$ is an {\em expansion} of $g_\theta$) i.e.
\begin{equation} d(x,y) \leq \rho (x,y), \;\;\; x,y \in M.
\label{e:30}
\end{equation}
(cf. \cite{kn:Str}, p. 230) where $d$ is the distance function
corresponding to the Webster metric. Although $\rho$ and $d$ are
inequivalent distance functions, they determine the same topology.
A first step towards recovering M. Obata's arguments (cf.
\cite{kn:Oba}) is the following
\begin{theorem}
\label{t:2} Let $(M , \theta )$ be a compact strictly pseudoconvex
CR manifold of CR dimension $n$, such that $\rho (X,X) + 2 A(X, J
X) \geq k G_\theta (X,X)$ for some $k
> 0$ and any $X \in H(M)$. Assume that
$\lambda \equiv - 2nk/(2n-1) \in {\rm Spec}(\Delta_b )$ and
${\mathcal H} \equiv {\rm Eigen}(\Delta_b ; \lambda ) \cap {\rm
Ker}(T) \neq (0)$. Then any eigenfunction $f \in {\mathcal H}$ is
given by
\begin{equation}
\label{e:31}
f(\gamma (s)) = \alpha \cos (s \sqrt{c}), \;\;\; s \in {\mathbb R},
\; c = k/(2n-1) \, ,
\end{equation}
along each lengthy geodesic $\gamma : {\mathbb R} \to M$ of the
Tanaka-Webster connection $\nabla$ such that $|\dot{\gamma}(s)| =
1$ and $\gamma (0) = x_0$, where $x_0 \in M$ is a point such that
$f(x_0 ) = \sup_{x \in M} f(x) \equiv \alpha$.

Assume additionally that $(M , \theta )$ is Sasakian {\rm (}$\tau
= 0${\rm )}. If any two points of $M$ can be joined by a
Carnot-Carath\'eodory minimizing lengthy geodesic then $f(x) =
\alpha \cos (r(x) \sqrt{c})$, $x \in M$, where $r(x) = \rho (x_0 ,
x)$ is the Carnot-Carath\'eodory distance from $x_0$. If $y_0 \in
M$ is a point such that $f(y_0 ) = \inf_{x \in M} f(x)$ then
$f(y_0 ) = - \alpha$. Consequently $M_{\pi /\sqrt{c}}$ consists
solely of critical points of $f$ and each $x \in M_{\pi
/\sqrt{c}}$ is degenerate. \label{t:final}
\end{theorem}
Here, for a given $s \in {\mathbb R}$ we let $M_s$ consist of all
points $x \in M$ such that there is a lengthy geodesic $\gamma :
{\mathbb R} \to M$ of $\nabla$, parametrized by arc length, such
that $\gamma (0) = x_0$ and $\gamma (s) = x$. The assumptions in
Theorem \ref{t:2} are rather restrictive and, among all odd
dimensional spheres, are satisfied only on $S^3$ (thus motivating
the conjecture in the Introduction). Precisely
\begin{proposition} Let $M = S^{2n+1}$ with the standard contact
form $\theta = \frac{i}{2} (\overline{\partial} - \partial )
|z|^2$. If {\rm i)} the inequality {\rm (\ref{e:8})} is satisfied
for some $k
> 0$, {\rm ii)} $- 2nk/(2n-1) \in {\rm Spec}(\Delta_b )$,  and {\rm iii)} ${\rm
Eigen}(\Delta_b ; - 2nk/(2n-1) ) \cap {\rm Ker}(T) \neq (0)$, then
$k = 4$ and $n=1$. Conversely the statements {\rm i)-iii)} hold on
$S^3$. Moreover if $M = S^3$ and
\[ \left. H = a(x_1^2 + y_1^2 -
x_2^2 - y_2^2 ) + 2 b(x_1 x_2 + y_1 y_2 ), \;\;\; f =
H\right|_{S^3} \, , \;\;\; b \neq 0, \] {\rm (}a spherical
harmonic of degree $2$ on $S^3$ such that $T(f) = 0${\rm )} then
$\alpha = \sup_{x \in S^3} f(x) = \sqrt{a^2 + b^2}$ and $f(\gamma
(s)) = \alpha \cos (2 s)$ for any lengthy geodesic $\gamma :
{\mathbb R} \to S^3$ of $\nabla$ {\rm (}the Tanaka-Webster
connection of $S^3${\rm )} parametrized by arc length and such
that $\gamma (0)$ is a maximum point of $f$. Moreover
\begin{equation}
M_{\pi /2} = \{ (\lambda , \mu ,  - \frac{b \lambda}{\alpha - a} ,
- \frac{b \mu}{\alpha - a}) : \lambda^2 + \mu^2 = \frac{\alpha -
a}{2 \alpha} \, , \;\; \lambda , \mu \in {\mathbb R} \}
\label{e:set}
\end{equation}
consists solely of degenerate critical points of $f$.
\label{p:nIS3}
\end{proposition}
The proof of Proposition \ref{p:nIS3} is relegated to Appendix A.
\par
{\em Proof of Theorem} \ref{t:final}. Assume that $\lambda = -
2nk/(2n-1)$ is an eigenvalue of $\Delta_b$ and let $f \in
{\mathcal H}$ be an eigenfunction of $\Delta_b$ corresponding to
$\lambda$ such that $T(f) = 0$. By the Bochner type formula
\eqref{e:7} one has
\[ \frac{1}{2} \, \Delta_b (|\nabla^H f |^2 ) = |\pi_H \nabla^2
f|^2 + \lambda |\nabla^H f|^2 + \rho (\nabla^H f \, , \, \nabla^H
f) + 2 L f \, .\] Once again we integrate and use Lemma \ref{l:2}
and the assumption \eqref{e:8}. We get
\[ 0 \geq \| \pi_H \nabla^2 f \|^2 - (1 + \frac{k}{\lambda}) \| \Delta_b f\|^2
= \| \pi_H \nabla^2 f \|^2 - \frac{1}{2n} \| \Delta_b f \|^2 \geq
0 \] (the last inequality is a consequence of \eqref{e:29}) hence
\[ \int_M \{ |\pi_H \nabla^2 f |^2 - \frac{1}{2n} (\Delta_b f)^2 \}
\Psi = 0 \] so that (again by \eqref{e:29})
\begin{equation}
\label{e:32}
|\pi_H \nabla^2 f|^2 = \frac{1}{2n} (\Delta_b f)^2 .
\end{equation}
The following lemma of linear algebra is well known. If $A \in
{\mathbb R}^{m^2}$ satisfies $m |A|^2 = {\rm trace}(A)^2$ then $A
= (1/m) \, {\rm trace}(A) I_m$, where $I_m$ is the unit matrix of
order $m$. Therefore (by \eqref{e:32})
\[ \pi_H \nabla^2 f = \frac{1}{2n} \, (\Delta_b f) \, G_\theta .
\]
In particular the identities \eqref{e:13}-\eqref{e:14} are
consistent with our assumption that $f_0 = 0$. Using again
$\Delta_b f = \lambda_1 f$ we may conclude that
\begin{equation}
\label{e:33}
\pi_H \nabla^2 f = - c \, f \, G_\theta \, ,
\end{equation}
where $c = k/(2n-1)$.
\par
M. Obata's proof (cf. {\em op. cit.}) of the fact that equality in
\eqref{e:1} yields $M^m \approx S^m$ (an isometry) is an
indication that we should evaluate \eqref{e:33} along a lengthy
geodesic of the Tanaka-Webster connection, and integrate the
resulting ODE. Let us recall briefly the needed material on
geodesics (as developed in \cite{kn:BaDr2}). Let $(U, x^1 , \cdots
, x^{2n+1})$ be a system of local coordinates on $M$ and let us
set $g \; d x^i = g^{ij}
\partial_j$, where $\partial_i = \partial /\partial x^i$. A
{\em sub-Riemannian geodesic} is a $C^1$ curve $\gamma (t)$ in
$M$ satisfying the Hamilton-Jacobi equations associated to the
Hamiltonian function $H(x, \xi ) = \frac{1}{2} \; g^{ij}(x) \xi_i
\xi_j$ that is
\begin{equation}
\label{e:34}
\frac{d x^i}{d t} = g^{ij}(\gamma (t)) \xi_j (t),
\end{equation}
\begin{equation}
\label{e:35}
\frac{d \xi_k}{d t} = - \frac{1}{2} \; \frac{\partial
g^{ij}}{\partial x^k}(\gamma (t)) \xi_i (t) \xi_j (t),
\end{equation}
for some cotangent lift $\xi (t) \in T^* (M)$ of $\gamma (t)$. Let
$\gamma (t) \in M$ be a sub-Riemannian geodesic and $s = \phi (t)$
a $C^1$ diffeomorphism. As shown in \cite{kn:BaDr2}, if $\gamma
(t) = \overline{\gamma}(\phi (t))$ then $\overline{\gamma}(s)$ is
a sub-Riemannian geodesic if and only if $\phi$ is affine, i.e.
$\phi (t) = \alpha t + \beta$, for some $\alpha , \beta \in
{\mathbb R}$. In particular, every sub-Riemannian geodesic may be
reparametrized by arc length $\phi (t) = \int_0^t
|\dot{\gamma}(u)| \, d u$. In \cite{kn:BaDr2} we introduced a {\em
canonical} cotangent lift of a given lengthy curve $\gamma : I \to
M$ by setting \[ \xi : I \to T^* (M), \;\;\; \xi (t) T_{\gamma
(t)} = 1, \;\;\; \xi (t) X = g_\theta (\dot{\gamma} , X), \] for
any $X \in H(M)_{\gamma (t)}$, and showed that
\begin{theorem}
\label{t:3} Let $M$ be a strictly pseudoconvex CR manifold and
$\theta$ a contact form on $M$ such that $G_\theta$ is positive
definite. A $C^1$ curve $\gamma (t) \in M$, $|t| < \epsilon$, is a
sub-Riemannian geodesic of $(M , H(M), G_\theta )$ if and only if
$\gamma (t)$ is a solution to
\begin{equation}
\label{e:36} \nabla_{\displaystyle{\dot{\gamma}}} \dot{\gamma} = -
2 b(t) J \dot{\gamma}, \;\;\; b^\prime (t) = A(\dot{\gamma},
\dot{\gamma}), \;\;\; |t| < \epsilon ,
\end{equation}
with $\dot{\gamma}(0) \in H(M)_{\gamma (0)}$, for some $C^1$
function $b : (-\epsilon , \epsilon ) \to {\mathbb R}$.
\end{theorem}
R. S. Strichartz's paper \cite{kn:Str} manifestly doesn't involve
any elements of connection theory or curvature. As argued by R. S.
Strichartz (cf. {\em op. cit.}) curvature is a measurement of the
deviation of the given Riemannian manifold from its Euclidean
model (and sub-Riemannian manifolds exhibit no approximate
Euclidean behavior). Nevertheless, in view of Theorem \ref{t:2}
when $(M, \theta )$ is a Sasakian manifold (i.e. $\tau = 0$) the
lengthy geodesics of $\nabla$ are among the sub-Riemannian
geodesics and it is likely that a variational theory of the
geodesics of $\nabla$ (as started in \cite{kn:BaDr2}) is the key
step towards bringing the results of \cite{kn:Oba} to CR geometry.
\par
Our approach (based on $\nabla$) is not in contradiction with the
arguments in \cite{kn:Str}: indeed the curvature of $\nabla$ is
related to the pseudoconvexity properties of $M$ (as understood in
complex analysis in several variables) rather than to its
intrinsic shape. To emphasize the impact of connection theory
within our approach we may prove the following elementary
regularity result. Note that a sub-Riemannian geodesic is required
to be of class $C^2$ (cf. \cite{kn:Str}, p. 233) and no higher
regularity is expected {\em a priori}. In turn, any $C^1$ geodesic
of $\nabla$ is automatically of class $C^\infty$ [as a projection
on $M$ of an integral curve of some standard horizontal vector
field (cf. Prop. 6.3 in \cite{kn:KoNo}, Vol. I, p. 139) having
$C^\infty$ coefficients].
\par
Let $\gamma (t) \in M$ be a lengthy geodesic of the Tanaka-Webster
connection, parametrized by arc-length ($|\dot{\gamma}(t)| = 1$).
Then (by \eqref{e:33})
\[ \frac{d^2 (f \circ \gamma )}{dt^2} = - c \; f \circ \gamma \]
hence $f(\gamma (t)) = A \cos (t \sqrt{c}) + B \sin (t \sqrt{c})$.
As $M$ is compact there is $x_0 \in M$ such that $f(x_0 ) =
\sup_{x \in M} f(x) =: \alpha$. Let $\gamma (t)$ be a lengthy
geodesic of $\nabla$ such that $\gamma (0) = x_0$. Then $A =
\alpha$ and $\{ d(f\circ \gamma )/d t \}_{t = 0} = 0$ yields $B =
0$ so that $f(\gamma (t)) = \alpha \; \cos (t \sqrt{c})$, which is
\eqref{e:31}.
\par
Again by compactness $(M , \rho )$ is a complete metric space,
hence (cf. Theorem 7.1 in \cite{kn:Str}, p. 244) any
sub-Riemannian geodesic can be extended indefinitely. Since $\tau
= 0$ the statements about sub-Rieman\-nian geodesics in
\cite{kn:Str} apply to the lengthy geodesics of $\nabla$ as well.
Let $\gamma : {\mathbb R} \to M$ be a lengthy geodesic of $\nabla$
such that $|\dot{\gamma}(s)| = 1$, $\gamma (0) = x_0$ and $\gamma
(s_{\rm min}) = y_0$. By \eqref{e:31}
\[ 0 = \frac{d}{ds} \{ f \circ \gamma \}_{s = s_{\rm min}} = -
\alpha \sqrt{c} \sin (\sqrt{c} \, s_{\rm min}) \] hence $s_{\rm
min} = m\pi /\sqrt{c}$ for some $m \in {\mathbb Z}$. Then $\alpha
> f(y_0 ) = (-1)^m \alpha$ implies that $m$ is odd.
Again by \eqref{e:31}, $M_{\pi /\sqrt{c}} \subset f^{-1} (-\alpha
)$. Finally, let $x \in M_{\pi/\sqrt{c}}$ and $(U, x^i )$ a local
coordinate system on $M$ such that $x \in U$. As $T(f) = 0$
\[ \frac{\partial^2 f}{\partial x^i \partial x^j}(x) T^i (x) T^j
(x) = (\nabla^2 f)(T,T)_x = T(f_0 )_x = 0 \] hence $x$ is a
degenerate critical point. Therefore, the points of $M_{\pi
/\sqrt{c}}$ may fail to be isolated. Nevertheless
\begin{proposition}
\label{p:1}
Let $(M , \theta )$ be a compact Sasakian
manifold. If for any $x \in B(x_0 , \pi /\sqrt{c})$ there is a
length minimizing {\rm (}with respect to the Carnot-Carath\'eodory
distance{\rm )} lengthy geodesic joining $x_0$ and $x$ then the
exponential map $\exp_{x_0} : N(x_0 , \pi /\sqrt{c}) \to B(x_0 ,
\pi /\sqrt{c})$ {\rm (}with respect to the Tanaka-Webster
connection{\rm )} is a surjection.
\end{proposition}
Here $B(x_0 , R) = \{ x \in M : \rho (x_0 , x) < R \}$ is the
Carnot-Carath\'eodory ball of center $x_0$ and radius $R > 0$.
Also $N(x_0 , R) = \{ w \in H(M)_{x_0} : |w| < R \}$.

{\em Proof of Proposition \ref{p:1}}. To see that the restriction
of $\exp_{x_0}$ to

\noindent $N(x_0 , \pi /\sqrt{c})$ is indeed $B(x_0 , \pi
/\sqrt{c})$-valued let $w \in N(x_0 , \pi /\sqrt{c})$, $w \neq 0$,
and $t = |w|$. Let us set $v = (1/t) w$ and consider the geodesic
$\gamma : {\mathbb R} \to M$ of $\nabla$ with the initial data
$\gamma (0) = x_0$ and $\dot{\gamma}(0) = v$, so that
\[ \exp_{x_0} (w) = \exp_{x_0} (t v) = \gamma (t). \]
Then
\[ \rho (x_0 , \gamma (t)) \leq \int_0^t |\dot{\gamma}(s)| ds = t
< \pi /\sqrt{c} \] i.e. $\gamma (t) \in B(x_0 , \pi /\sqrt{c})$.
\par
To see that $\exp_{x_0} : N(x_0 , \pi /\sqrt{c}) \to B(x_0 , \pi
/\sqrt{c})$ is on-to let $x \in B(x_0 , \pi /\sqrt{c})$ and let
$\gamma : {\mathbb R} \to M$ be a length minimizing geodesic
joining $x_0$ and $x$ and such that $\dot{\gamma}(0) = v \in
H(M)_{x_0}$, with $|v| = 1$. Then $\gamma (t) = x$ for some $t \in
{\mathbb R} \setminus \{ 0 \}$, so that $\exp_{x_0} (t v) = x$.
Next
\[ g_{\theta , x_0} (t v , t v) = t^2 |v|^2 = r(x)^2 < \pi^2 /c \]
i.e. $t v \in N(x_0 , \pi /\sqrt{c})$. Q.e.d.
\par
A generalization of the Lichnerowicz-Obata theorem (\cite{kn:Lic},
\cite{kn:Oba}) to the case of Riemannian foliations was obtained
by J.M. Lee \& K. Richardson, \cite{kn:LeRi} (see also
\cite{kn:LeRi0}). The leaf space of a Riemannian foliation is
often an orbifold (for instance if all leaves are compact) so that
(in light of \cite{kn:DrMa}) one expects analogs to Theorems
\ref{t:1} and \ref{t:2} on a CR orbifold (see also E. Stanhope,
\cite{kn:Sta}). This matter will be addressed in a further paper.
\begin{appendix}
\section{On the spectrum of the sublaplacian on the standard sphere} Let
$M$ be a strictly pseudoconvex CR manifold and $\theta$ a contact
form on $M$ with $G_\theta$ positive definite. Let $\nabla^\theta$
be the Levi-Civita connection of the semi-Riemannian manifold $(M,
g_\theta )$. Then (cf. e.g. \cite{kn:DrTo}, Chapter 1)
\begin{equation}
\nabla^\theta = \nabla + (\Omega - A) \otimes T + \tau \otimes
\theta + 2 \theta \odot J \label{e:A.1}
\end{equation}
where $\odot$ is the symmetric tensor product. Then
\begin{equation} \label{e:A.5}
\nabla^\theta_X X = \nabla_X X - A(X,X) T, \;\;\; X \in H(M).
\end{equation}
Given a local $G_\theta$-orthonormal frame $\{ X_a : 1 \leq a \leq
2n \}$ of $H(M)$ one has (by (\ref{e:A.5}) and ${\rm trace}(\tau )
= 0$)
\[ \Delta f = \sum_{j = 0}^{2n} \{ X_j (X_j f) -
(\nabla^\theta_{X_j} X_j ) (f) \} = T(T(f)) + \Delta_b f \] for
any $f \in C^\infty (M)$, where $X_0 = T$, proving Greenleaf's
formula (\ref{e:Greenleaf}). Let ${\mathcal P}_\ell$ be the set of
all homogeneous polynomials $H : {\mathbb R}^{2n+2} \to {\mathbb
R}$ of degree ${\rm deg}(H) = \ell$ and ${\mathcal H}_\ell =
{\mathcal P}_\ell \cap {\rm Ker}(\Delta_{{\mathbb R}^{2n+2}})$. To
compute eigenvalues of $\Delta_b$ starting from ${\rm Spec}(\Delta
)$ we consider the equation
\begin{equation}
\Delta_b f + T^2 (f) = - \ell (2n + \ell ) f \label{e:A.2}
\end{equation}
with $f = \left. H\right|_{S^{2n+1}}$ and $H \in {\mathcal
H}_\ell$. For example if $\ell = 1$ and $H \in {\mathcal H}_1 =
{\mathcal P}_1$ then $T_0^2 (H) = - H$ hence $-2n \in {\rm
Spec}(\Delta_b )$. In general
\begin{proposition} If there is $\lambda \in {\mathbb R}$ and ${\mathcal H}_\ell \cap
{\rm Ker}(T_0^2 + \lambda I) \neq (0)$ then $\lambda - \ell (2n +
\ell ) \in {\rm Spec}(\Delta_b )$. For instance one may produce
the eigenvalues $\{ - 2n , \; - 4n , \; - 6n - 8, \; - 6n \}
\subset {\rm Spec}(\Delta_b )$ and ${\rm Eigen}(\Delta_b ; - 2n)
\cap {\mathcal P}_1 \neq (0)$, ${\rm Eigen}(\Delta_b ; - 4n) \cap
{\mathcal P}_2 \neq (0)$ and ${\rm Eigen}(\Delta_b ; \lambda )
\cap {\mathcal P}_3 \neq (0)$ for each $\lambda \in \{ - 6n - 8,
\; - 6 n \}$. \label{p:2}
\end{proposition}
If $H = \sum_{i,j=1}^{n+1} (a_{ij} x^i x^j + b_{ij} x^i y^j +
c_{ij} y^i y^j ) \in {\mathcal P}_2$ (with $a_{ij}, b_{ij} \in
{\mathbb R}$, $a_{ji} = a_{ij}$, $c_{ji} = c_{ij}$) then $T_0^2 H
= - \lambda H$ if and only if $2(c_{ij} - a_{ij}) = - \lambda
a_{ij}$, $2(b_{ij} + b_{ji}) = \lambda b_{ij}$ and $2(c_{ij} -
a_{ij}) = \lambda c_{ij}$. Hence ${\rm Ker}(T_0^2 + \lambda I)
\cap {\mathcal P}_2 = (0)$ for any $\lambda \in {\mathbb R}
\setminus \{ 4 \}$ and
\[ {\rm Ker}(T_0^2 + 4 I) \cap {\mathcal P}_2 = \] \[ = \{ a_{ij} (x^i x^j
- y^i y^j ) + b_{ij} x^i y^j : a_{ij} , b_{ij} \in {\mathbb R}, \;
a_{ij} = a_{ji} \} \subset {\mathcal H}_2 . \] Similarly ${\rm
Ker}(T_0^2 + \lambda I) \cap {\mathcal P}_3 = (0)$ for any
$\lambda \in {\mathbb R} \setminus \{ 1, 9 \}$ and
\[ {\rm Ker}(T_0^2 + I) \cap {\mathcal H}_3 = \{  (a_{ijk} x^i + b_{ijk} y^i )(x^j x^k + y^j y^k )
: \] \[ a_{ijk}, b_{ijk} \in {\mathbb R} \;\; {\rm symmetric}, \;
\sum_j a_{ijj} = \sum_j b_{ijj} = 0, \; 1 \leq i \leq n+1 \} , \]
\[ {\rm Ker}(T_0^2 + 9 I) \cap {\mathcal H}_3 = \{
a_{ijk} x^i (x^j x^k - 3 y^j y^k ) + b_{ijk} (y^i y^j - 3 x^i x^j
) y^k :
\]
\[ a_{ijk}, b_{ijk} \in {\mathbb R} \;\; {\rm symmetric}, \;
\sum_j a_{ijj} = \sum_j b_{ijj} = 0, \; 1 \leq i \leq n+1 \} .
\]
Proposition \ref{p:2} is proved. The calculation of the full ${\rm
Spec}(\Delta_b )$ on $S^{2n+1}$ is an open problem. \vskip 0.1in
{\em Proof of Proposition} \ref{p:nIS3}. Let $R$ be the curvature
of the Tanaka-Webster connection. Then (cf. Chapter 1 in
\cite{kn:DrTo})
\begin{equation}
R(X,Y)Z = G_\theta (Y,Z) X - G_\theta (X,Z) Y + \label{e:curv}
\end{equation}
\[ + G_\theta (J Y, Z) J X - G_\theta (J X , Z) J Y - 2 G_\theta
(J X , Y ) J Z \] for any $X,Y,Z \in H(S^{2n+1})$. Taking the
trace in (\ref{e:curv}) we obtain
\begin{equation}
\rho (X,X) = 2(n+1) G_\theta (X,X). \label{e:Ricci}
\end{equation}
The assumptions i)-ii) imply that $-2nk/(2n-1)$ is an eigenvalue
of the ordinary Laplacian on $S^{2n+1}$. On the other hand
${\mathcal H}_1 \cap {\rm Ker}(T_0 ) = (0)$ hence $2nk/(2n-1)$ is
greater equal than $4(n+1)$. Finally (by (\ref{e:Ricci})) $k \leq
2(n+1)$ hence $n = 1$ and $k = 4$.
\par
Let $\nabla^\theta$ be the Levi-Civita connection of $S^3$. As
$S^3$ is a Sasakian manifold $\nabla^\theta_{\dot{\gamma}}
\dot{\gamma} = \nabla_{\dot{\gamma}} \dot{\gamma} + 2 \theta
(\dot{\gamma}) J \dot{\gamma}$ (by (\ref{e:A.1})) for any $C^1$
curve $\gamma (t)$ in $S^3$. In particular any lengthy geodesic
$\gamma$ of $\nabla$ is a geodesic of $S^3$ as well. Moreover any
geodesic $\gamma$ of $\nabla$ with $\dot{\gamma}(0) \in H(S^3
)_{\gamma (0)}$ is lengthy. Indeed (as $\nabla T = 0$)
\[ \frac{d}{d t} \{ \theta_{\gamma (t)} (\dot{\gamma}(t)) \} =
g_\theta (\dot{\gamma} , \nabla_{\dot{\gamma}} T)_{\gamma (t)} = 0
\]
hence $\theta (\dot{\gamma})_{\gamma (t)} = \theta_{\gamma (0)}
(\dot{\gamma}(0)) = 0$. Let $x_0 \in S^3$ such that $\alpha =
f(x_0 )$. Let $\gamma$ be a lengthy geodesic of $\nabla$,
parametrized by arc length, such that $\gamma (0) = x_0$. Then
$\gamma (s) = x_0 \cos s + x \sin s$, $s \in {\mathbb R}$, for
some $x \in {\mathbb R}^4$ such that $\| x \| = 1$ and $\langle
x_0 , x \rangle = 0$. If $U = S^3 \setminus \{ x_2 = y_2 = 0 \}$
the Levi distribution $\left. H(S^3 )\right|_U$ is spanned by
\[ X = \frac{\partial}{\partial x_1} - F \,
\frac{\partial}{\partial x_2} - G \, \frac{\partial}{\partial y_2}
\, , \;\;\; Y = \frac{\partial}{\partial y_1} + G \,
\frac{\partial}{\partial x_2} - F \, \frac{\partial}{\partial y_2}
\, , \]
\[ F(x,y) = \frac{x_1 x_2 + y_1 y_2}{x_2^2 + y_2^2} \, , \;\;\;
G(x,y) = \frac{x_1 y_2 - y_1 x_2}{x_2^2 + y_2^2} \, , \] hence the
condition that $\gamma$ is lengthy reads
\[ \left. Q^j \frac{\partial}{\partial x^j}\right|_{x_0} + R^j
\left. \frac{\partial}{\partial y^j} \right|_{x_0} = \lambda
X_{x_0} + \mu Y_{x_0} \] for some $\lambda , \mu \in {\mathbb R}$,
where $Q^j = x^j (x)$ and $R^j = y^j (x)$, $j \in \{ 1,2 \}$, or
\begin{equation}
Q^1 = \lambda , \;\;\; Q^2 = \mu G(x_0 ) - \lambda F(x_0 ) ,
\label{e:43}
\end{equation}
\begin{equation}
R^1 = \mu , \;\;\; R^2 = - \mu F(x_0 ) - \lambda G(x_0 ).
\label{e:44}
\end{equation}
Let us set $P^j = x^j (x_0 )$ and $S^j = y^j (x_0 )$. The solution
to the constrained extreme value problem $\alpha = \sup_{x \in
S^3} f(x)$ is $\alpha = \sqrt{a^2 + b^2}$ and
\[ P^1 = \xi , \;\; S^2 = \eta , \;\; P^2 = A \xi , \;\; S^2 = A
\eta , \;\;\; \xi^2 + \eta^2 = \frac{\alpha -a}{2\alpha} ,
\] where $A = (\alpha - a)/b$, hence $F(x_0 ) = b/(\alpha - a)$
and $G(x_0 ) = 0$. Finally $\| x \| = 1$ may be written $\lambda^2
+ \mu^2 = (\alpha - a)/(2\alpha )$ hence (\ref{e:43})-(\ref{e:44})
yield (\ref{e:set}) in Proposition \ref{p:nIS3}.
\end{appendix}

\vskip 0.5cm {\small Universit\`a degli Studi della Basilicata,

Dipartimento di Matematica,

Campus Macchia Romana,

85100 Potenza, Italy,

e-mail: {\tt barletta@unibas.it}}

\begin{thebibliography}{12345}
\bibitem [1]{kn:BaDr} Barletta, E. and Dragomir, S., On the
spectrum of a strictly pseudoconvex CR manifold, Abhandlungen
Math. Sem. Univ. Hamburg, {\bf 67}(1997), 143-153.
\bibitem [2]{kn:BaDr2} Barletta, E. and Dragomir, S., Jacobi fields
of the Tanaka-Webster connection on Sasakian manifolds, to appear
in Kodai Mathematical Journal, 2006.
\bibitem [3]{kn:BGM} Berger, M.,  Gauduchon, P. and Mazet, E.,
Le spectre d'une variet\'e Riemannienne, Lecture Notes in Math.,
Vol.  {\bf 194}, Springer-Verlag, Berlin-New York, 1971.
\bibitem [4]{kn:Bla} Blair, D. E., Contact manifolds in Riemannian
geometry, Lecute Notes in Math., Vol. {\bf 509}, Springer-Verlag,
Berlin-Heidelberg-New York, 1976.
\bibitem [5]{kn:Boc} Bochner, S., Vector fields and Ricci
curvature, Bull. Amer. Math. Soc., {\bf 52}(1946), 776-797.
\bibitem [6]{kn:Cho} Chow, W. L.,  \"Uber Systeme von Lineaaren
Partiellen Differentialgleichungen erster Ordnung, Math. Ann. {\bf
117} (1939), 98-105.
\bibitem [7]{kn:Dra} Dragomir, S.,  A survey of pseudohermitian
geometry, The Proceedings of the Workshop on Differential Geometry
and Topology, Palermo (Italy), June 3-9, 1996, in Supplemento ai
Rendiconti del Circolo Matematico di Palermo II, {\bf 49}(1997),
101-112.
\bibitem [8]{kn:DrUr} Dragomir, S. and Urakawa, H., On the inhomogeneous
Yang-Mills equation $d^*_D R^D = f$, Interdisciplinary Information
Sciences, (1){\bf 6}(2000), 41-52.
\bibitem [9]{kn:DrMa} Dragomir, S. and Masamune, J.,
Cauchy-Riemann orbifolds, Tsukuba J. Math., (2){\bf 26}(2002),
351-386.
\bibitem [10]{kn:DrNi} Dragomir, S. and Nishikawa, S., Foliated CR
manifolds, J. Math. Soc. Japan, (4){\bf 56}(2004), 1031-1068.
\bibitem [11]{kn:DrTo} Dragomir, S.  and Tomassini, G.,
Differential geometry and analysis on CR manifolds, Progress in
Mathematics, Vol. {\bf 246}, Birkh\"auser, Boston-Basel-Berlin,
2006.
\bibitem [12]{kn:Gre} Greenleaf, A., The first eigenvalue of a
sublaplacian on a pseudohermitian manifold, Commun. Partial
Differential Equations, {\bf 10}(1985), 191-217.
\bibitem [13]{kn:JiHo} Jiaqing, Z. and Hongcang, Y.,  On the estimate
of the first eigenvalue of a compact Riemannian manifold, Scientia
Sinica, (12){\bf 27}(1984), 1265-1273.
\bibitem [14]{kn:KoNo} Kobayashi, S. and Nomizu, K.,
Foundations of differential geometry, Interscience Publishers, New
York, Vol. I, 1963, Vol. II, 1969.
\bibitem [15]{kn:Lee} Lee, J. M.,  Pseudo-Einstein structures on
CR manifolds, Amer. J. Math., {\bf 110}(1988), 157-178.
\bibitem [16]{kn:LeRi0} Lee, J. M. and Richardson, K., Riemannian
foliations and eigenvalue comparison, Annals of Global Analysis
and Geometry, {\bf 16}(1998), 497-525.
\bibitem [17]{kn:LeRi} Lee, J. M., and Richardson, K., Lichnerowicz
and Obata theorems for foliations, Pacific Journal of Mathematics,
(2){\bf 206}(2002), 339-357.
\bibitem [18]{kn:Lic} Lichnerowicz, A., G\'eom\'etrie des
groupes de transformations, Dunod, Paris, 1958.
\bibitem [19]{kn:MaYe} Marinescu, G. and Yeganefar, N.,
Embeddability of some strongly pseudoconvex CR manifolds,
preprint, {\tt arXiv:math.CV/0403044 v1}, {\tt 2 Mar 2004}.
\bibitem [20]{kn:MeSj} Menikoff, A. and Sj\"ostrand, J.,  On the
eigenvalues of a class of hypoelliptic operators, Math. Ann., {\bf
235}(1978), 55-58.
\bibitem [21]{kn:Oba} Obata, M., Certain conditions for a
Riemannian manifold to be isometric with a sphere, J. Math. Soc.
Japan, {\bf 14}(1962), 333-340.
\bibitem [22]{kn:Sta} Stanhope, E., Spectral bounds on orbifold
isotropy, Annals of Global Analysis and Geometry, {\bf 27}(2005),
355-375.
\bibitem [23]{kn:Str} Strichartz, R. S., Sub-Riemannian
geometry, J. Differential Geometry, {\bf 24}(1986), 221-263.
\bibitem [24]{kn:Tan} Tanaka, N., A differential geometric study
on strongly pseudo-convex manifolds, Kinokuniya Book Store Co.,
Ltd., Kyoto, 1975.
\bibitem [25]{kn:Web} Webster, S. M., Pseudohermitian structures
on a real hypersurface, J. Differential Geometry, {\bf 13}(1978),
25-41.
\end{thebibliography}
\end{document}